\documentclass{article}

\begin{document}

\begin{center}
{\LARGE Symmetry Groups and Equivalence Transformations\smallskip }

{\LARGE in the Nonlinear Donnell--Mushtari--Vlasov Theory\smallskip }

{\LARGE for Shallow Shells}{\Large \footnote{%
Published in \textit{Journal of Theoretical and Applied Mechanics, 1997,
Year XXVII, No. 2, 43-51. }}}{\LARGE \bigskip }

{\Large Vassil M.Vassilev}

Institute of Mechanics, Bulgarian Academy of Sciences

Acad. G. Bontchev St., Block 4, 1113 Sofia, Bulgaria

E-mail: vassil@bgcict.acad.bg{\Large \medskip }
\end{center}

\begin{quotation}
\textbf{Abstract--}In the case of transversely only loaded shallow shells,
the nonlinear Donnell--Mushtari--Vlasov theory for large deflection of
isotropic thin elastic shells leads to a system of two coupled nonlinear
forth-order partial differential equations known as Marguerre's equations.
This system involves two arbitrary elements -- the curvature tensor of the
shell middle-surface and the function of transversal load per unit surface
area. In the present note, the point symmetry groups of Marguerre's
equations are established, the corresponding group classification problem
being solved. It is shown that Marguerre's equations are equivalent to the
von K\'{a}rm\'{a}n equations for large deflection of plates in the
time-independent case and in the time-dependent case as well. It is also
observed that the same holds true in respect of the field equations for
anisotropic shallow shells.\smallskip 
\end{quotation}

\section{Introduction}

Within the framework of the nonlinear Donnell--Mushtari--Vlasov (DMV) theory
(see, e.g., \cite{NIORD}) the state of equilibrium of a transversely loaded
thin isotropic elastic shell of uniform thickness is determined by the
following system of two coupled nonlinear forth-order partial differential
equations: 
\begin{equation}
\begin{array}{l}
D\Delta ^{2}w-\varepsilon ^{\alpha \mu }\varepsilon ^{\beta \nu }w_{;\alpha
\beta }\Phi _{;\mu \nu }-\varepsilon ^{\alpha \mu }\varepsilon ^{\beta \nu
}b_{\alpha \beta }\Phi _{;\mu \nu }=p, \\ 
(1/Eh)\Delta ^{2}\Phi +(1/2)\varepsilon ^{\alpha \mu }\varepsilon ^{\beta
\nu }w_{;\alpha \beta }w_{;\mu \nu }+\varepsilon ^{\alpha \mu }\varepsilon
^{\beta \nu }b_{\alpha \beta }w_{;\mu \nu }=0,
\end{array}
\label{EQ1}
\end{equation}
in two independent variables, associated with the coordinates on the shell
middle-surface $F,$ and two dependent variables -- the transversal
displacement function $w,$ and Airy's stress function $\Phi $. Here, $%
\varepsilon ^{\alpha \beta }$ is the alternating tensor of $F$; $b_{\alpha
\beta }$ is the curvature tensor of $F$; $D$, $E$ and $h$ are the bending
rigidity, Young's modulus and thickness of the shell, respectively (i.e., $D$%
, $E$ and $h$ are given constants); $p$ is the function of transversal load
per unit surface area; a semicolon is used for covariant differentiation
with respect to the metric tensor $a_{\alpha \beta }$ of the surface $F$; $%
\Delta $ is the Laplace-Beltrami operator on $F$. Here and throughout, Greek
indices have the range 1, 2 and the usual summation convention over a
repeated index (one subscript and one superscript) is employed.

The present note is concerned with a special case of the nonlinear DMV
theory; the so-called \textit{shallow shells} are considered. In fact, this
is an approximation of the theory, which is usually introduced as follows
(cf., e.g. \cite{GALIM}). Let $(x^1,x^2,x^3)$ be a fixed right-handed
rectangular Cartesian coordinate system in the 3-dimensional Euclidean space
in which the middle-surface $F$ of a shell is embedded, and let this surface
be given by the equation 
\[
x^3=f(x^1,x^2),\;(x^1,x^2)\in \Omega \subset \mathbf{R}^2, 
\]
where $f:\mathbf{R}^2\rightarrow \mathbf{R}$ is assumed to be a
single-valued and smooth function possessing as many derivatives as may be
required on the domain $\Omega $. Let us take $x^1,x^2$ to serve as
coordinates on the surface $F$. Then, relative to this coordinate system,
the components of the fundamental tensors and the alternating tensor of $F$
are given by the expressions: 
\begin{equation}
a_{\alpha \beta }=\delta _{\alpha \beta }+f_{,\alpha }f_{,\beta },\quad
b_{\alpha \beta }=a^{-1/2}f_{,\alpha \beta },\quad \varepsilon ^{\alpha
\beta }=a^{-1/2}e^{\alpha \beta },  \label{EQ2}
\end{equation}
where 
\[
a=\det (a_{\alpha \beta })=1+(f_{,1})^2+(f_{,2})^2; 
\]
$\delta _{\alpha \beta }=\delta ^{\alpha \beta }$ is the Kronecker delta
symbol; $e^{\alpha \beta }$ is the alternating symbol; here and in what
follows, a comma is used for partial differentiation with respect to the
coordinates on\ $F$. \textit{A shell is said to be shallow on the domain} $%
\Omega _0\subset \Omega $\ \textit{when the inequalities} 
\[
\left| f_{,\alpha }\right| \left| f_{,\beta }\right| \leq \varepsilon ^2\ll
1,\quad \varepsilon =const, 
\]
\textit{hold for every point} $(x^1,x^2)\in \Omega _0$. Hence, for shallow
shells the quadratic terms in the right-hand sides of expressions (2) are
small compared to unity and may be neglected. Thus, allowing for a relative
error of order $O\left( \varepsilon ^2\right) $, one may regard the
intrinsic geometry of the shell middle-surface $F$ as Euclidean and $%
(x^1,x^2)$ may be thought of as an Euclidean coordinate system in which: 
\begin{eqnarray}
a_{\alpha \beta } &=&\delta _{\alpha \beta },  \label{EQ3a} \\
b_{\alpha \beta } &=&f_{,\alpha \beta },  \label{EQ3b} \\
\varepsilon ^{\alpha \beta } &=&e^{\alpha \beta };  \label{EQ3c}
\end{eqnarray}
the mean curvature $H$ of the surface $F$ and its Gaussian curvature $K$
(note that the latter is not necessarily equal to zero within the allowed
relative error) take the form 
\begin{eqnarray}
H &=&(1/2)\delta ^{\alpha \beta }f_{,\alpha \beta },  \label{EQ4a} \\
K &=&(1/2)e^{\alpha \mu }e^{\beta \nu }f_{,\alpha \beta }\/f_{,\mu \nu },
\label{EQ4b}
\end{eqnarray}
and system (\ref{EQ1}) reads 
\begin{equation}
\begin{array}{l}
D\delta ^{\alpha \beta }\delta ^{\mu \nu }w_{,\alpha \beta \mu \nu
}-e^{\alpha \mu }e^{\beta \nu }w_{,\alpha \beta }\Phi _{,\mu \nu }-e^{\alpha
\mu }e^{\beta \nu }b_{\alpha \beta }\Phi _{,\mu \nu }=p, \\ 
(1/Eh)\delta ^{\alpha \beta }\delta ^{\mu \nu }\Phi _{,\alpha \beta \mu \nu
}+(1/2)e^{\alpha \mu }e^{\beta \nu }w_{,\alpha \beta }w_{,\mu \nu
}+e^{\alpha \mu }e^{\beta \nu }b_{\alpha \beta }w_{,\mu \nu }=0,
\end{array}
\label{EQ5}
\end{equation}
Thus, we arrive at the \textit{equilibrium equations for shallow shells}
within\ the\ framework\ of\ the\ nonlinear\ DMV\ theory.\ Since\ equations\ (%
\ref{EQ5})\ follow\ from\ Marguerre's\ shell\ theory\ \cite{MARG}\ as\
well,\ they\ are\ also\ known\ as \textit{Marguerre's equations} for large
deflection of plates with small initial curvature (i.e., shallow shells).
These equations are well accepted and play an important role in the shell
theory (see, e.g., \cite{YYYU}, \cite{NIKH} and the references therein).
They also include as a special case, with $b_{\alpha \beta }=0,$ the
well-known \textit{von K\'{a}rm\'{a}n equations} for large deflection of
plates \cite{KARM}.

\section{Symmetry Groups}

The aim of the present work is to study, following \cite{OVS}, \cite{OLV}
and \cite{IBR}, the invariance properties of system (\ref{EQ5}) relative to
local one-parameter Lie groups of local point transformations acting on open
subsets of the 4-dimensional Euclidean space $\mathbf{R}^4$, with
coordinates $(x^1,x^2,w,\Phi )$, representing the involved independent and
dependent variables. For that purpose Lie infinitesimal technique is used
and, as a rule, the results obtained are expressed in terms of the \textit{%
infinitesimal generators (operators)} of the groups; in the present case,
the latter are vector fields of the form 
\begin{equation}
\mathbf{X}=\xi ^{\mu } \frac \partial {\partial x^{\mu } }+\eta \frac
\partial {\partial w}+\varphi \frac \partial {\partial \Phi },  \label{EQ6}
\end{equation}
where $\xi ^{\mu } ,\eta $ and $\varphi $ are functions of the variables $%
x^1,x^2,w$ and $\Phi .$ The system considered involves an arbitrary tensor
field -- $b_{\alpha \beta }$, and an arbitrary function -- $p$. This gives
rise to a \textit{group classification problem} with respect to the \textit{%
arbitrary element} -- the set $\left\{ b_{\alpha \beta },p\right\} $.

The infinitesimal criterion of invariance leads to the following system of
determining equations (\textit{DE system}) for the components $\xi ^1,\xi
^2,\eta $ and $\varphi $ of the vector fields of form (\ref{EQ6}) generating
point symmetry groups admitted by system (\ref{EQ5}): 
\begin{equation}
\xi ^1=C_1x^1+C_2x^2+C_3,  \label{EQ7}
\end{equation}
\begin{equation}
\xi ^2=-C_2x^1+C_1x^2+C_4,  \label{EQ8}
\end{equation}
\begin{equation}
\eta =\eta (x^1,x^2),  \label{EQ9}
\end{equation}
\begin{equation}
\varphi =B_1x^1+B_2x^2+B_3,  \label{EQ10}
\end{equation}
\begin{equation}
b_{\alpha \mu }\xi _{,\beta }^{\mu } +b_{\beta \mu }\xi _{,\alpha }^{\mu }
+\xi ^{\mu } b_{\alpha \beta ,\mu }=-\eta _{,\alpha \beta },  \label{EQ11}
\end{equation}
\begin{equation}
e^{\alpha \mu }e^{\beta \nu }b_{\alpha \beta }\eta _{,\mu \nu }=0,
\label{EQ12}
\end{equation}
\begin{equation}
D\delta ^{\alpha \beta }\delta ^{\mu \nu }\eta _{,\alpha \beta \mu \nu
}=2p\xi _{,\mu }^{\mu } +\xi ^{\mu } p_{,\mu },  \label{EQ13}
\end{equation}
where $C_1,...,C_4,B_1,B_2,$ and $B_3$ are arbitrary real constants. This
result is established in \cite{VAS} through the standard computational
procedure (see \cite[Sec. 5]{OVS} or \cite[Sec. 2.4]{OLV}), and therefore we
can assert that the \textit{system (\ref{EQ5}) admits a vector field of form
(\ref{EQ6}) if and only if (\ref{EQ7}) -- (\ref{EQ13}) hold}. In other
words, all possible symmetries of the kind under consideration inherent to
system (\ref{EQ5}) can be found via the solution of the DE system. So, the
problem to solve consists in finding the solutions of this system. The main
difficulty here comes out of the determining equations (\ref{EQ11}) -- (\ref
{EQ13}). They show that the space of solutions $L$ to the DE system, i.e.,
the Lie algebra of symmetries associated with system (\ref{EQ5}), depends on
the choice of the arbitrary element $\{b_{\alpha \beta },p\}$. At this
juncture, we do face a group classification problem. That is to say, that we
will have to determine all those \textit{specializations} (special forms) of
the arbitrary element for which system (\ref{EQ5}) admits vector fields of
the form (\ref{EQ6}).

Proceeding to analyze this problem, we will first show that equations (\ref
{EQ9}), (\ref{EQ11}) -- (\ref{EQ13}) may be replaced by the following three
equivalent ones: 
\begin{equation}
\eta =-\xi ^{\mu }f_{,\mu }+A_1x^1+A_2x^2+A_3,  \label{EQ14}
\end{equation}
where $A_1,\,A_2$ and $A_3$ $\ $are arbitrary real constants, and 
\begin{equation}
2P\xi _{,\mu }^{\mu } +\xi ^{\mu } P_{,\mu }=0,  \label{EQ15}
\end{equation}
\begin{equation}
2K\xi _{,\mu }^{\mu } +\xi ^{\mu } K_{,\mu }=0,  \label{EQ16}
\end{equation}
where 
\begin{equation}
P=2D\delta ^{\mu \nu }H_{,\mu \nu }+p.  \label{EQ17}
\end{equation}
Indeed, (\ref{EQ14}) represents the general solution of equations (\ref{EQ11}%
) of the form (\ref{EQ9}) when (\ref{EQ3b}), (\ref{EQ7}) and (\ref{EQ8})
hold; under the same assumption, by substituting (\ref{EQ14}) into (\ref
{EQ12}) and (\ref{EQ13}), and taking into account (\ref{EQ4a}) and (\ref
{EQ4b}), after some algebra we obtain expressions (\ref{EQ15}) and (\ref
{EQ16}). Using this result, hereafter we will assume that \textit{the DE
system consists of equations (\ref{EQ7}), (\ref{EQ8}), (\ref{EQ10}), (\ref
{EQ14}), (\ref{EQ15}) and (\ref{EQ16})}.

Now, a 6-dimensional space of solutions to the DE system arises immediately.
It corresponds to the 6-dimensional Lie algebra $L_0$ of the vector field
whose components are given as follows 
\begin{eqnarray}
\xi ^{\mu } &=&0,  \nonumber  \label{Eq18} \\
\eta &=&A_1x^1+A_2x^2+A_3,  \label{EQ18} \\
\varphi &=&B_1x^1+B_2x^2+B_3,  \nonumber
\end{eqnarray}
where $A_1,\,A_2,\,A_3,\,B_1,\,B_2$ and $B_3$ are arbitrary real constants.
Evidently, the above solutions of the DE system do not depend on the choice
(specialization) of the arbitrary element. Hence, \textit{for any
specialization of the arbitrary element system (\ref{EQ5}) admits the
6-parameter Lie group }$G_0$\textit{\ generated by all linear combinations
of the vector fields} 
\[
\frac \partial {\partial w},\quad x^1\frac \partial {\partial w},\quad
x^2\frac \partial {\partial w},\quad \frac \partial {\partial \Phi },\quad
x^1\frac \partial {\partial \Phi },\quad x^2\frac \partial {\partial \Phi }; 
\]
these vector fields constitute a basis of the associated Lie algebra $L_0$.
Simultaneously, $G_0$ \textit{is the largest Lie group of point
transformations admitted by system (\ref{EQ5}) for any choice of the
arbitrary element}. Indeed, by setting 
\[
f=\sin (x^1)\sin (x^2), 
\]
for instance, one can easily verify that for this particular choice of the
arbitrary element the DE system has no other solutions besides those given
by (\ref{EQ18}).

So far, we have obtained the so-called \textit{kernel of the full symmetry
groups} associated with system (\ref{EQ5}), that is the group $G_0$. The
next step is to identify the cases in which the system under consideration
possesses larger groups of point symmetries. In the light of all the above,
it means to characterize in a suitable manner all those specializations of
the arbitrary element for which system (\ref{EQ5}) admits vector fields of
the form 
\begin{equation}
\mathbf{Y}=\xi ^{\mu } \frac \partial {\partial x^{\mu } }-\xi ^{\mu }
f_{,\mu }\frac \partial {\partial w}\quad (\xi ^{\mu } \neq 0),  \label{EQ19}
\end{equation}
where $\xi ^{\mu } $ are given by the expressions (\ref{EQ7}) and (\ref{EQ8}%
), keeping in mind that \textit{(\ref{EQ15}) and (\ref{EQ16}) remain the
only conditions (necessary and sufficient) for (\ref{EQ5}) to be invariant
under a group generated by a vector field of form (\ref{EQ19})}. Note that
all determining equations are thus taken into account.

Taking into account (\ref{EQ4a}), (\ref{EQ4b}) and (\ref{EQ17}) we can see
at once that for $f=p=0$ the invariance conditions (\ref{EQ15}) and (\ref
{EQ16}) are satisfied. Hence, system (\ref{EQ5}), with $f=p=0$, possesses a
larger group of point symmetries (in addition to $G_0$) and (\ref{EQ7}), (%
\ref{EQ8}) and (\ref{EQ19}) show that this is the complete 4-parameter group
of homothetic motions of the Euclidean plane. In this special case, (\ref
{EQ5}) coincide with the homogeneous von K\'arm\'an equations, so that we
have arrived at the result obtained in \cite{SHW} (see also \cite{AMES}).
This example gives us a good motivation for studying the general case.

We begin with the following observation. \textit{Given a vector field} $%
\mathbf{Y}${}\textit{of form (\ref{EQ19}), the function} 
\begin{equation}
\tilde w=w+f,  \label{EQ20}
\end{equation}
\textit{is an invariant of the corresponding Lie group of transformations
acting on} $\mathbf{R}^4$. Therefore, we can introduce new coordinates $%
(x^1,x^2,\tilde w,\Phi )$ on $\mathbf{R}^4$, the new dependent variable $%
\tilde w$ being defined by (\ref{EQ20}), in which $\mathbf{Y}$\textbf{\ }%
takes the form of an infinitesimal operator of the Lie algebra associated
with the group of homothetic motions of the Euclidean plane, namely 
\begin{equation}
\mathbf{Y}=\xi ^{\mu } \frac \partial {\partial x^{\mu } },  \label{EQ21}
\end{equation}
(note that the components $\xi ^1$ and {$\xi ^2$} still have the form (\ref
{EQ7}) and (\ref{EQ8}), respectively, as before), and system (\ref{EQ5})
reads 
\begin{equation}
\begin{array}{l}
D\delta ^{\alpha \beta }\delta ^{\mu \nu }\tilde w_{,\alpha \beta \mu \nu
}-e^{\alpha \mu }e^{\beta \nu }\tilde w_{,\alpha \beta }\Phi _{,\mu \nu }=P,
\\ 
(1/Eh)\delta ^{\alpha \beta }\delta ^{\mu \nu }\Phi _{,\alpha \beta \mu \nu
}+(1/2)e^{\alpha \mu }e^{\beta \nu }\tilde w_{,\alpha \beta }\tilde w_{,\mu
\nu }=K.
\end{array}
\label{EQ22}
\end{equation}
Thus, \textit{the problem of invariance of system (\ref{EQ5}) with respect
to a vector field (\ref{EQ19}) converts into the problem of invariance of
system (\ref{EQ22}) under a vector field (\ref{EQ21})} as a change of the
variables does not affect the group properties of a system of differential
equations. Now, taking into account the invariance conditions (\ref{EQ15})
and (\ref{EQ16}) (note that they remain unchanged under the above coordinate
transformation), we can conclude that \textit{system (\ref{EQ22}) admits a
one-parameter group }$G$\textit{\ of homothetic motions of the Euclidean
plane generated by a vector field of form (\ref{EQ21}) if and only if the
corresponding arbitrary element }$\{b_{\alpha \beta },p\}$\textit{\ is such
that the functions }$P$\textit{\ and }$K$\textit{, defined by formulae (\ref
{EQ4b}) and (\ref{EQ17}), are invariants of }$G$\textit{\ (when }$C_1=0$%
\textit{) or eigenfunctions (when }$C_1\neq {0}$\textit{) of its generator,
the latter being regarded as an operator acting on the smooth functions} $%
\zeta :M\rightarrow R,\;M\subset \mathbf{R}^2$. This result may be thought
of as a general solution to the group classification problem under
consideration in terms of the function $p$ and the two characteristic
invariants, $H$ and $K$, of the shell middle-surface $F$.

\section{Equivalence Transformations}

Let us now briefly discuss, in the context of the shell theory, the meaning
of the coordinate transformation $\omega :\mathbf{R}^4\rightarrow \mathbf{R}%
^4,\;(x^1,x^2,w,\Phi )\longmapsto (x^1,x^2,w+f,\Phi )$ introduced in the
previous Section. Omitting tilde's in equations (\ref{EQ22}), we can say
that systems (\ref{EQ5}) and (\ref{EQ22}) belong to the same class -- they
have the same differential structure, and differ from one another only in
the form of the arbitrary element. Hence, according to \cite[Definition 6.4]
{OVS}, we may conclude that $\omega :\mathbf{R}^4\rightarrow \mathbf{R}^4$ 
\textit{is an equivalence transformation for system(\ref{EQ5})}. In
addition, we see that (\ref{EQ22}) is nothing but a system of nonhomogeneous
von K\'arm\'an equations with special right-hand sides. It is noteworthy
that \textit{Marguerre's equations (\ref{EQ5}), i.e., the equilibrium
equations for shallow shells, turned out to be equivalent to the von
K\'arm\'an equations, i.e., the equilibrium equations for plates},
regardless of the invariance properties of system (\ref{EQ5}). To the best
of the author's knowledge, this fact has not been noticed before in the
literature, though both the von K\'arm\'an equations and Marguerre's
equations have been studied and utilized for many years until now (see,
e.g., \cite{NIORD}, \cite{GALIM}, \cite{YYYU}, \cite{NIKH}, \cite{SHW}, \cite
{AMES}, \cite{CIA} and the references therein). In our opinion, the
established correspondence between these two systems of equations will
certainly be of use in solving a wide range of problems arising in the shell
theory. In particular, the results for the von K\'arm\'an equations obtained
in \cite{DJO} can be easily conveyed to the theory of shallow shells.

Several additional applications of the transformation $\omega $ can by given
too. Evidently, it maps:

$\bullet $ the time-dependent Marguerre's equations 
\[
\begin{array}{l}
D\delta ^{\alpha \beta }\delta ^{\mu \nu }w_{,\alpha \beta \mu \nu
}-e^{\alpha \mu }e^{\beta \nu }w_{,\alpha \beta }\Phi _{,\mu \nu }-e^{\alpha
\mu }e^{\beta \nu }b_{\alpha \beta }\Phi _{,\mu \nu }+\rho \ddot w=p, \\ 
(1/Eh)\delta ^{\alpha \beta }\delta ^{\mu \nu }\Phi _{,\alpha \beta \mu \nu
}+(1/2)e^{\alpha \mu }e^{\beta \nu }w_{,\alpha \beta }w_{,\mu \nu
}+e^{\alpha \mu }e^{\beta \nu }b_{\alpha \beta }w_{,\mu \nu }=0,
\end{array}
\]
(where $\rho $ is the mass per unit area of the shell middle-surface, and a
superposed dot is used to denote partial derivative with respect to the time 
$t$) into the time-dependent von K\'arm\'an equations 
\[
\begin{array}{l}
D\delta ^{\alpha \beta }\delta ^{\mu \nu }w_{,\alpha \beta \mu \nu
}-e^{\alpha \mu }e^{\beta \nu }w_{,\alpha \beta }\Phi _{,\mu \nu }+\rho
\ddot w=P, \\ 
(1/Eh)\delta ^{\alpha \beta }\delta ^{\mu \nu }\Phi _{,\alpha \beta \mu \nu
}+(1/2)e^{\alpha \mu }e^{\beta \nu }w_{,\alpha \beta }w_{,\mu \nu }=K;
\end{array}
\]

$\bullet $ the time-dependent equations for anisotropic shallow shells 
\[
\begin{array}{l}
D^{\alpha \beta \mu \nu }w_{,\alpha \beta \mu \nu }-e^{\alpha \mu }e^{\beta
\nu }w_{,\alpha \beta }\Phi _{,\mu \nu }-e^{\alpha \mu }e^{\beta \nu
}b_{\alpha \beta }\Phi _{,\mu \nu }+\rho \ddot w=p, \\ 
E^{\alpha \beta \mu \nu }\Phi _{,\alpha \beta \mu \nu }+(1/2)e^{\alpha \mu
}e^{\beta \nu }w_{,\alpha \beta }w_{,\mu \nu }+e^{\alpha \mu }e^{\beta \nu
}b_{\alpha \beta }w_{,\mu \nu }=0,
\end{array}
\]
(where $D^{\alpha \beta \mu \nu }$ and {$E^{\alpha \beta \mu \nu }$ }denote
the material constants) into the following system of von K\'arm\'an-type
equations 
\[
\begin{array}{l}
D^{\alpha \beta \mu \nu }w_{,\alpha \beta \mu \nu }-e^{\alpha \mu }e^{\beta
\nu }w_{,\alpha \beta }\Phi _{,\mu \nu }+\rho \ddot w=D^{\alpha \beta \mu
\nu }f_{,\alpha \beta \mu \nu }+p, \\ 
E^{\alpha \beta \mu \nu }\Phi _{,\alpha \beta \mu \nu }+(1/2)e^{\alpha \mu
}e^{\beta \nu }w_{,\alpha \beta }w_{,\mu \nu }=K.
\end{array}
\]

The above list could be continued.\bigskip

\noindent \textbf{\Large Acknowledgments}\medskip

This research was supported by Contract No. MM 517/1995 with the National
Scientific Fund, Republic of Bulgaria.\newpage\


\begin{thebibliography}{99}
\bibitem{NIORD}  \textsc{Niordson, F. I.} \textit{Shell Theory.}
Nord-Holland, Amsterdam (1985).

\bibitem{GALIM}  \textsc{Galimov, K. Z.} \textit{Foundations of the
Nonlinear Theory of Shells.} Kazan' University Press, Kazan (1975).

\bibitem{MARG}  \textsc{Marguerre, K.} Zur Theorie der gekr\"{u}mmten Platte
gro\ss er Form\"{a}nderung. Proc. Fifth Intern. Congr. Appl. Mech., 93,
Cambridge, Massachusetts (1938).

\bibitem{YYYU}  \textsc{Yu, Y. Y.} On equations for large deflections of
elastic plates and shallow shells. \textit{Mech. Res. Commun.}, \textbf{18},
(1991) pp. 373-384.

\bibitem{NIKH}  \textsc{Nakazawa, M., T. Iwakuma, S. Kuranishi and M. Hudaka.%
} Instability phenomena of a rectangular elastic plate under bending and
shear. \textit{Int. J. Solids Structures}, \textbf{30} (1993), pp. 2729-2741

\bibitem{KARM}  \textsc{von Karman, Th.} Festigkeitesprobleme im
Maschinebau, Encyklop\"{a}die der Mathematischen Wissenschaften. Bd IV, 311,
Taubner, Leipzig (1910).

\bibitem{OVS}  \textsc{Ovsiannikov, L. V.} \textit{Group Analysis of
Differential Equations.} Nauka, Moscow (1978). English transl. edited by W.
F. Ames, Academic Press, New York (1982).

\bibitem{OLV}  \textsc{Olver, P. J.} \textit{Applications of Lie Groups to
Differential Equations.} 2nd edition. Graduate Texts in Mathematics, Vol.
107, Springer-Verlag, New York (1993).

\bibitem{IBR}  \textsc{Ibragimov, N. Kh.} \textit{Transformation Groups
Applied to Mathematical Physics}. Nauka, Moscow (1983). English transl.,
Riedel, Boston (1985).

\bibitem{VAS}  \textsc{Vassilev, V. M.} Group Analysis of a Class of
Equations of the Theory of Thin Elastic Plates and Shells. Ph.D. thesis,
Bulgarian Academy of Sciences (1991).

\bibitem{SHW}  \textsc{Shwarz, F.} Lie symmetries of the von K\'{a}rm\'{a}n
equations, \textit{Comput. Phys. Commun.}, \textbf{31 }(1984), pp. 113-114.

\bibitem{AMES}  \textsc{Ames, K. A., W. F. Ames.} On group analysis of the
von K\'{a}rm\'{a}n equations. \textit{Nonlinear Analysis}, \textbf{6 }%
(1982), pp. 854-853.

\bibitem{CIA}  \textsc{Ciarlet, Ph., P. Rabier.} \textit{Les Equations de
von K\'{a}rm\'{a}n, }Lecture Notes in Math., Vol. 826, Springer-Verlag,
Berlin-Heidelberg-New York (1980).

\bibitem{DJO}  \textsc{Djondjorov, P., V. Vassilev.} Conservation laws and
group-invariant solutions of the von K\'{a}rm\'{a}n equations. \textit{Int.
J. Non-Linear Mechanics}, \textbf{31} (1996), pp. 73-87.
\end{thebibliography}
\end{document}